\documentclass[leqno, 12pt, twoside]{article}

\usepackage{amsmath}
\usepackage{amssymb}
\usepackage{amsthm}  
\usepackage{latexsym}

\newcommand{\E}{{\cal E}}
\newcommand{\D}{{\cal D}}
\newcommand{\F}{{\cal F}}
\newcommand{\A}{{\cal A}}

\newcommand{\G}{{\cal G}}

\newcommand{\s}{\mathrm{span}}
\newcommand{\aff}{\mathrm{aff}}

\newcommand{\bu}{{\mathbf u}}\newcommand{\bv}{{\mathbf v}}
\newcommand{\bw}{{\mathbf w}}\newcommand{\bx}{{\mathbf x}}
\newcommand{\by}{{\mathbf y}}

\newcommand{\n}{{\mathbf n}}

\newcommand{\CC}{{\mathbf C}}
\newcommand{\NN}{{\mathbf N}}
\newcommand{\RR}{{\mathbf R}}

\DeclareMathOperator{\codim}{codim}

\DeclareMathOperator{\sign}{sign}

\DeclareMathOperator{\Hom}{Hom}
\DeclareMathOperator{\affspan}{span}

\newcommand{\mm}{{\mathfrak{m}}}

\newcommand{\id}{\mathrm{id}}

\newcommand{\into}{\hookrightarrow}
\newcommand{\longto}{\longrightarrow}
\newcommand{\qr}{\overline}
\newcommand{\quer}[1]{\,\overline{\!#1}} 
\newcommand{\ba}{\mathbf{a}}

\newcommand{\x}{\times}

\renewcommand{\:}{\colon}

\renewcommand{\epsilon}{\varepsilon} 
\renewcommand{\phi}{\varphi}
\renewcommand{\rho}{\varrho}
\renewcommand{\theta}{\vartheta}

\newtheorem{theorem}{Theorem}[section]
\newtheorem{corollary}[theorem]{Corollary}
\newtheorem{lemma}[theorem]{Lemma}
\newtheorem{definition}[theorem]{Definition}
\newtheorem{remark}[theorem]{Remark}
\newtheorem{remdef}[theorem]{Remark and Definition}
\newtheorem{proposition}[theorem]{Proposition}
\newtheorem{HRsetup}[theorem]{Abstract HR setup}

\newtheorem*{HLT}{Hard Lefschetz Theorem (HLT)} 
\newtheorem*{HRR}{Hodge-Riemann Bilinear Relations (HRR)} 

\title{Hodge-Riemann Relations for Polytopes\\
A Geometric Approach\\}
\author{Gottfried Barthel, Jean-Paul Brasselet,\\
Karl-Heinz Fieseler, Ludger Kaup}
\date{}
\begin{document}
\maketitle
\begin{abstract}
    \noindent
The key to the Hard Lefschetz Theorem for combinatorial
intersection cohomology of polytopes is to prove the Hodge-Riemann
bilinear relations. In these notes, we strive to present an easily
accessible proof. The strategy essentially follows the original
approach of \cite{karu}, applying induction \textit{\`a la}
\cite{brel2}, but our guiding principle here is to emphasize the
geometry behind the algebraic arguments by consequently stressing
polytopes rather than fans endowed with a strictly convex conewise
linear function. It is our belief that this approach makes the
exposition more transparent since polytopes are more appealing to
our geometric intuition than convex functions on a~fan.
\end{abstract}

\section{Introduction}

The proof of the Hard Lefschetz Theorem for the
``Combinatorial Intersection Cohomo\-logy'' of
polytopes given in \cite{karu} was the keystone
in a long endeavour of several research groups to
verify that Stanley's generalized (``toric'')
$h$-vector for polytopes has the conjectured
properties: The theorem (usually referred to as
``{\bf HLT}'' in the sequel) implies that the
generalized $h$-vector agrees with the vector of
even degree Intersection Cohomology Betti numbers
and that this vector enjoys the unimodality property
(in addition to symmetry and non-negativity).

The HLT is an easy consequence of the so-called
bilinear ``Hodge-Riemann relations'' (``HR
relations'' or ``{\bf HRR}'' for short); and
since the latter, being a ``positivity result",
reflect convexity in a more appropriate way than
the HLT, the focus has shifted towards proving
these relations. The first proof of the HRR
given in \cite{karu} has been rather involved.
The task of making it more easily accessible has
been taken up in different articles, cf.\
\cite{brel2} and \cite{mars}. With the present
notes, we further pursue this direction: Being
convinced that polytopes are closer to our
geometric intuition, we present an approach that
stresses geometric operations on polytopes
rather than algebraic operations on strictly
convex conewise linear functions.

Let us briefly recall the setup, referring to section~4 for
further details: To an $n$-dimensional polytope~$P$ in an
$n$-dimensional real vector space~$V$, one associates its
outer normal fan $\Delta = \Delta(P)$ in the dual vector
space $V^*$, and a conewise linear strictly convex
function~$\psi$. The ``combinatorial intersection cohomology''
$IH(\Delta)$ is a finite-dimensional real vector space with
even grading \( \bigoplus_{k = 0}^n IH^{2k}(\Delta) \). There
is a perfect pairing
$$
\cap \: IH^q(\Delta) \times IH^{2n-q}(\Delta) \longto \RR\ ,
$$
the ``intersection product'', so Poincar\'e duality holds on
$IH(\Delta)$.

On \( IH(\Delta) \), the multiplication with~$\psi$ induces
an endomorphism
\[
  L : IH^q(\Delta) \longto IH^{q+2}(\Delta)
\]
called the \textbf{Lefschetz operator}. The key result
of \cite{karu} (see also \cite{brel2}) reads as follows:
\begin{HLT}
    For each $k \geqq 0$, the iterated Lefschetz operator
\[
  L^k \: IH^{n-k}(\Delta) \longto IH^{n+k}(\Delta)
\]
is an isomorphism.
\end{HLT}

\noindent By Poincar\'e duality, it suffices to prove that each
map $L^k$ be injective or surjective.

\medskip
Using the intersection product, the Hard Lefschetz Theorem
can be restated in a different framework: Each mapping~$L^k$
(for $k \geqq 0$) yields a bilinear form
\[
  s_{k} \:  IH^{n-k}(\Delta) \x IH^{n-k}(\Delta) \longto \RR
  \;,\quad (\xi,\eta) \longmapsto \xi \cap L^k(\eta),
\]
called the $k$-th \textbf{Hodge-Riemann bilinear form}, or
\textbf{``HR-form''} for short. This form is symmetric since~$L$ is
self-adjoint with respect to the intersection product. In this
set-up, the HLT is equivalent to the non-degeneracy of all
forms~$s_{k}$.

Beyond non-degeneracy, the HR relations provide explicit
formul{\ae} for the signatures of these pairings. To that
end, we have to consider the \textbf{primitive intersection
cohomology}
\[
  IP^{n-k}(\Delta) \; : = \; \ker\bigl(L^{k+1}:{IH^{n-k}(\Delta)}
  \longto {IH^{n+k+2}(\Delta)}\bigr)
\]
for \( 0 \leqq k \leqq n \) (with \( k \equiv n \mod 2 \)).
In fact, assuming the HLT, there
is an $s_k$-orthogonal decomposition
$$
  IH^{n-k}(\Delta) \;=\; L \bigl(IH^{n-k-2}(\Delta)\bigr)
  \oplus IP^{n-k}(\Delta)\,.
$$
More generally, we see:

\begin{proposition}
\label{orthosum} If the HLT holds for the Lefschetz operator $L$
on the intersection cohomology of the fan $\Delta = \Delta(P)$,
then, for each $k$, the intersection cohomology splits as an
orthogonal direct sum
\[
  IH^{n-k}(\Delta) \;=\; \bigoplus_{j \geqq 0} \,
  L^j\bigl(IP^{n-k-2j}(\Delta)\bigr) \,.
\]
\end{proposition}

Now for each $q \leqq n-2$, the restricted operator~\( L \)
provides an isometric embedding ${IH^q(\Delta)} \into IH^{q+2}(\Delta)$
with respect to the pertinent HR-forms. Hence, in order to
determine the signature of $s_k$, it
suffices to consider the restrictions of the Hodge Riemann
forms $s_{k+2j}$ to the corresponding primitive subspaces
$IP^{n-k-2j}(\Delta)$. Here is the statement:

\begin{HRR}
    For each $k \geqq 0$ (with \( k \equiv n \mod2 \)), the
    Hodge-Riemann bilinear form~$s_{k}$ is
    $(-1)^{(n-k)/2}$-definite on~$IP^{n-k}(\Delta)$\,.
\end{HRR}

The HR relations imply the HLT, since the HR forms are readily
seen to be non-degenerate by descending induction on~$k$: For
$k = n$, that follows from $IP^0(\Delta) = IH^0(\Delta)$. For
\( k < n \), we assume that $s_{k+2}$ is non-degenerate. Then
so is the restriction
of $s_k$ to $L\bigl(IH^{n-k-2}(\Delta)\bigr)$.
This implies
$$
  IH^{n-k}(\Delta) = L(IH^{n-k-2}(\Delta))
  \oplus L\bigl(IH^{n-k-2}(\Delta)\bigr)^\perp,
$$
and it now suffices to prove
$L\bigl(IH^{n-k-2}(\Delta)\bigr)^\perp =
IP^{n-k}(\Delta)$. The inclusion ``$\supset$"
follows from the fact that $L$ is
$\cap$-self-adjoint, while ``$\subset$" is a
consequence of Poincar\'e duality for the
complementary dimensions $n-k-2$ and $n+k+2$.
\medskip

From Proposition \ref{orthosum}, we immediately
obtain a reformulation of the HRR in which the
primitive cohomology does not enter explicitly:

\begin{proposition}
\label{HRequat} The HRR are equivalent to the
HLT together with the additional condition that
the Hodge-Riemann bilinear forms~$s_{k}$ on
$IH^{n-k}(\Delta)$ satisfy the \textrm{``{\bf
HR-equation}''}
    \[
      \sign(s_{k}) = \sign(s_{k+2}) + (-1)^{(n-k)/2}(b_{n-k} -
      b_{n-k-2}) \, ,
    \]
    where $b_q :=  \dim_{\RR} IH^q(\Delta)$ denotes the $q^{{\rm th}}$
intersection cohomology Betti number of the fan
$\Delta$.
\end{proposition}

\section{Outline of the proof of the Hodge-Riemann relations}
The HR relations are known to hold if the polytope $P$ is simple.
The first proof has been given in \cite{mcmu}; a simplified
version followed in \cite{timo}. This result is the basis for
the proof of the general case by a twofold induction: The
``outer loop''  is on the dimension $n :=  \dim(P)$. For the
more involved ``inner loop'', following \cite{brel2}, we
associate to~$P$ an integer $\mu :=  \mu(P) \geqq 0$ that
measures how far~$P$ is from being simple: It counts
those faces, here called ``normally stout'' (see
\ref{stout}), that witness non-simplicity, with $\mu = 0$
characterizing simple polytopes. The inner induction
on~$\mu$ requires three main steps:

\begin{description}
    \item[\textsc{Cutting off}] (Section \ref{S3}):
    Given a face $F \prec P$, we consider an affine
    hyperplane $H$ that is sufficiently near and parallel
    to a supporting hyperplane for the face $F \prec P$
    and intersects $\,\overset{\circ}{\!P}$. Let $P =
    G \cup R$ be the corresponding decomposition of $P$ into
    the ``germ~$G = G_P (F)$ of $P$ along the face $F$" and
    the residual polytope~$R$.

    Then, if $F \prec P$ is a normally stout face of minimal
    dimension, we have $\mu(R) < \mu(P)$, so the HRR hold
    for~$R$ by induction hypothesis.
        For the investigation of the germ $G$, it is
        important that the face $F$ itself is a simple
        polytope and that it is ``normally trivial"
    in~$P$, cf.\ \ref{normtrivmin}.

    \item[\textsc{HRR for special $\boldsymbol{n}$-polytopes}]
    (Section \ref{Perprop}):
    Assuming the HRR to hold for $m$-poly\-to\-pes (with $m <
    n$), we prove the validity for the following special
    $n$-polytopes:
     \begin{itemize}
       \item[5.1] A pyramid~$P = \Pi(Q)$ with an $(n-1)$-dimensional
       base $Q$.

    \item[5.2] A non-trivial product $P = S \x
        P_{0}$, where~$S$ is simple.

    \item[] Furthermore, we prove the following \textbf{Gluing
    property}:

    \item[5.3]  The HRR hold for an $n$-polytope~$P$ that can
    be cut ``transversally'' into two polytopes~$P_{1}$
    and~$P_{2}$ such that the HRR hold for both pieces.
    \end{itemize}

    \item[\textsc{Deformation}]
    of the germ $G$ into a product (Section \ref{Deformation}):
    There is a continuous family $(Q_{t})_{t \in [0,1]}$ of
    pairwise combinatorially equivalent polytopes with $Q_1 = G$
    and $Q_0 = F \x {\Pi}(L)$ with the pyramid~$  \Pi(L)$ over a
    ``link"~$L = L_P (F)$ of $F$ in $P$.
    Then HRR is valid for~$Q_1 = G$  iff it is for~$Q_{0}$.

\end{description}

By induction hypothesis, the HRR hold for the
lower
    dimensional polytopes $F$ and~$L$, and thus, by 5.1
    and 5.2, hold also for $Q_0 = F \x {\Pi}(L)$, hence
    eventually also for $G$.
Finally, the gluing result of 5.3 applied with
$P_1 = R$ and $P_2 = G$ from Step~1 (``Cutting
off'') yields the HRR for the initial
polytope~$P$.

\medskip
In section \ref{intersect}, we recall the
definition and basic properties of combinatorial
intersection cohomology as needed later on.

\section{Cutting off}
\label{S3} In this section, we explain how a
polytope can be made simple by successively
cutting off faces containing non-simple points.
In that process, we have to make sure at each
step that we get closer to the class of simple
polytopes. A measure for the ``distance" of a
polytope~\( P \) to that class is the number
$\mu (P)$ of its ``normally stout" faces, see
Def.~\ref{stout}.

\smallskip
We first introduce some basic constructions.

\begin{remdef}
\label{cutting} The complement $V \setminus H$
of an affine hyperplane $H \subset V$ consists
of two open connected components
$$
  V \setminus H = U_1 \cup U_2\, .
$$
We say that a subset $A \subset V$ \emph{lies
strictly on one side} of $H$ if either $A
\subset U_1$ or $A \subset U_2$.

    Let \( P \subset V \) be a polytope, and~\( H \), a
    hyperplane as above. We call~\( H \) a \emph{cutting
    hyperplane} for~\( P \) if it intersects the relative
    interior, i.e., $H \cap \overset{\circ}{P} \not =
    \emptyset$. Such a hyperplane yields a
    decomposition
$$
  P = P_1 \cup P_2
$$
of~\( P \) into polytopes $P_i := P \cap
\quer{U}\!_i$ with \( \dim P_{1} = \dim P_{2} =
\dim P \). Both pieces meet along the common
facet \( P_{1} \cap P_{2} = H \cap P \) which we
call the \emph{cut facet}.

We say that \( H \) cuts \( P \)
\emph{transversally} if no vertex of~\( P \)
lies on~\( H \). Moreover, for a proper face \(
F \precneqq P \), we say that such a transversal
hyperplane~\( H \) is \emph{sufficiently near}
to~\( F \) (or a ``nearby hyperplane'') if~\( F
\) lies strictly on one side of~\( H \), whereas
all the remaining vertices of~\( P \) (i.e.,
those not contained in~\( F \)) lie on the other
side. If in addition $H$ is parallel to a
supportimg hyperplane $H_0 \subset V$ for the
face $F \prec P$, i.e. $P \cap H_0 = F$, we also
say that~\( H \) ``cuts off'' the face~\( F \).
\end{remdef}

In the sequel, cutting off a proper face \( F \)
of~\( P \) by a nearby parallel hyperplane plays
an important role: The resulting decomposition
\[
P=G \cup R
\]
of $P$ into one polytope $G$ containing the
face~\( F \) and a ``residual polytope'' $R$
allows a ``\emph{divide et impera}'' approach to
the HRR problem.

\begin{definition}
\label{germlink} Let $F \precneqq P$ be a proper
face of a polytope $P$ in $V$.
\begin{enumerate}
\item A {\bf germ} $G = G_P(F)$ of $P$ along the
face $F$ is any polytope $G$ obtained as
follows: Choose an affine hyperplane $H_0$ in
$V$ with $P \cap H_0 = F$ and let $H$ be a
parallel hyperplane such that $F$ lies strictly
on one side, say $U_1$, of $H$, and the vertices
of $P$ not contained in $F$ on the other side
$U_2$. Then $G:=  P_1$, while $R:= P_2$ is the
{\bf corresponding residual polytope}.

\item A {\bf link} $L = L_P(F)$ of the face $F
\precneqq P$ is any polytope obtained in the
following way:
\begin{itemize}
\item If $F = \{ {\ba} \}$, then $L_P(F) :=
L_P({\ba}) := P \cap H$ (the cut facet) for $H$
as above.

\item If $\dim F > 0$, choose a transversal
affine subspace $N \subset V$ to~\( F \), i.e.,
complementary to the affine span $\aff (F)$, and
intersecting the relative interior
$\overset{\circ}{F}$. Put $\widetilde P := P
\cap N$ and $\widetilde F := F \cap N$. Then
$L_P(F) := L_{\widetilde P}( \widetilde F)$ (so
\( L_P(F) = N \cap P \cap H \)).
\end{itemize}
\end{enumerate}
\end{definition}

We note that the combinatorial type of a germ
and that of a link is independent of all choices
made in the construction. --- In the literature
on convex polytopes, a link of a vertex is
usually called ``\emph{vertex figure}'', and a
link \( L_{P}(F) \) of a face is called
``\emph{face figure}'' or ``\emph{quotient
polytope}'', often noted \( P/F \).

\smallskip
We recall the notion of the \emph{join} of two
``relatively skew'' polytopes.

\begin{definition}
\label{join} Let $Q_{1}, Q_{2}$ be disjoint
polytopes in $V$ such that $\dim \aff (Q_{1}
\cup Q_{2}) = \dim Q_{1} + \dim Q_{2} + 1$. Then
their {\bf join} $Q_{1} * Q_{2}$ is the convex
hull of $Q_{1} \cup Q_{2}$ in~\( V \).
\end{definition}

We note that \( Q * \emptyset = Q \), and \(
Q_{1} * Q_{2} = Q_{2} * Q_{1} \). We remark that
the join $Q_{1} * Q_{2}$ is the disjoint union
of $Q_{1}$, $Q_{2}$, and all open segments
$(x,y)$ joining points $x \in Q_{1}$ and $y \in
Q_{2}$. We further remark that all faces of the
join are of the form \( F_{1} * F_{2} \), where
\( F_{i} \preceq Q_{i} \) is a (possibly empty)
face, and that a link \( L_{Q_{1} *
Q_{2}}(Q_{1}) \) is combinatorially equivalent
to~\( Q_{2} \).
--- We denote with $\Pi (P) := P * \{{\ba}\}$ for
${\ba} \not\in \aff (P)$ the {\bf pyramid} with
apex ${\ba}$ and base~$P$. An iterated pyramid
$\Pi^ i (P)$ for $i> 0$ is thus a join $P *
S_{i-1}$ with an $(i\,{-}\,1)$-simplex
$S_{i-1}$, whereas $\Pi^0(P) = P$.

\smallskip
We now study the local geometry near a face $F$
of a polytope $P$ in~$V$. For a given vertex
${\ba} \in F$, we fix a nearby cutting
hyperplane $H \subset V$. The cut facet $H \cap
P$ is a link $L_P ({\ba})$ of ${\ba}$ in~$P$,
and its face $F \cap H$ is a link $L_F({\ba})$
of ${\ba}$ relative to~\( F \).

\begin{definition}
\label{normtrivial} A proper face $F \precneqq
P$ of the polytope $P$ is called \vspace{-6pt}
\begin{itemize}
\item {\bf normally trivial (in \( P \)) at the
vertex \,${\ba}$\,} if the link $L_P({\ba})$ is
the join $L_F({\ba}) * S_{\ba}$ with a suitable
``complementary'' face $S_{\ba} \preceq
L_P({\ba})$\,, and \vspace{-6pt} \item {\bf
normally trivial (in $P$)} if it is normally
trivial at each of its vertices.
\end{itemize}
\end{definition}

We remark that every vertex of a polytope~\( P
\) is normally trivial as a face. If~\( {\ba} \)
is a simple vertex of~\( P \), then every
face~\( F \) containing~\( {\ba} \) is normally
trivial at~\( \ba \): A link \( L = L_{P}(\ba)
\) of~\( \ba \) in~\( P \) is a simplex, so for
the face \( F' = L_{F}(\ba) \) of~\( L \), there
is a unique complementary face; the latter being
again a simplex, any link of~\( F \) in~\( P \)
is a simplex. If \( F \precneqq P \) is an edge
or a facet of a three-dimensional polytope, then
the converse holds: Normal triviality at a
vertex \( {\ba} \) is equivalent to \( {\ba} \)
being simple.

More generally, for a face \( F \precneqq P \)
that is normally trivial at the vertex~\( {\ba}
\), there is a unique face  \( F_{\ba}'
\precneqq P \) ``cutting out'' the complementary
face \( S_{\ba} \) in the link \( L \), i.e.,
satisfying \( S_{\ba} = L_{P}({\ba}) \cap
F'_{\ba}\). That face is complementary to~\( F
\) at~\( {\ba} \), i.e., we have \( F \cap
F_{\ba}' = {\ba} \), \( \dim F + \dim F_{\ba}' =
\dim P \), and every edge emanating from \(
{\ba} \) either lies in~\( F \) or in \(
F_{\ba}' \). Shifting the affine span of \(
F_{\ba}' \) to the relative interior of~\( F \)
yields a transversal subspace~\( N \) to~\( F \)
as in Def.\ \ref{germlink},~2. As a consequence,
the polytope $S_{\ba}$ has the same
combinatorial type as $L_P(F)$, so that type
does not depend on the vertex ${\ba} \in F$.

\smallskip
Normal triviality of a face yields a
combinatorial local product structure:

\begin{remark}
\label{combiprod} Let $F \precneqq P$ be a
normally trivial face with link $L = L_P(F) = N
\cap P \cap H$ as in Def.\ \ref{germlink},~2.
Denote with $G = G_P(F)$ a corresponding germ
and  with $\pi \: V \to N$ the (affine)
projection onto $N$ along $\aff (F)$, i.e.,
collapsing \( \aff (F) \) to a single point~\(
\bv_0 \). Then $\pi$ induces a surjective map
$$
  \pi|_G \: G \longto \Pi(L)
$$
onto the pyramid $\Pi(L) := G \cap N$ over $L$
with apex $\bv_0$, mapping vertices onto
vertices. Moreover, we obtain a bijection
between the vertices of $G$ and the vertices of
$F \times \Pi(L)$ as follows: A vertex $\bu$
lying on the ``ridge'' \( F  \) of the ``hip
roof" $G$ is mapped to $(\bu,\bv_0)$, and a
vertex $\bv$ lying on the ``bottom facet'' \( G
\cap H \) (i.e., on the cut facet), being the
end point of an edge emanating from a unique
vertex $\bu \in F$, is mapped to $(\bu, \pi
(\bv))$. That map yields a \emph{combinatorial
equivalence} between the polytopes $G$ and $F
\times \Pi(L)$.
\end{remark}

If the link of a face is not a pyramid, then no
vertex lying on that face is a simple point of
the ambient polytope. This observation motivates
the interest of the following concept that is
essential for the inner loop, cf. \cite{brel2}
2.7:

\begin{definition}
\label{stout}
\begin{enumerate}

\item
  A (non-empty) polytope $P$ is called\/ {\bf stout}
  if it is not the pyramid over one of its facets.

\item A face $F \preceq P$ is called\/
  {\bf normally stout in $\boldsymbol{P}$} if one
  (and thus any) link $L = L_P(F)$ is stout.
\end{enumerate}
\end{definition}

Equivalently, a polytope~\( P \) is stout if for
each facet~\( F \), there are at least two
vertices of~\( P \) not lying on~\( F \). Hence,
``stoutness" only depends on the combinatorial
type and $\dim P \ge 2$ for a stout
polytope~$P$. In particular, a normally stout
face \( F \prec P \) always has codimension at
least $3$.

\smallskip
The relation general versus stout polytopes is
as follows:

\begin{lemma}\label{stoutfactor}
    If a polytope~\( P \) is not a simplex, then it has
    exactly one maximal stout face~\( B \preceq
    P\).
    In particular~\( P \) is the
    iterated pyramid
    \[
      P = \Pi^c(B) = B * S_{c-1}
    \]
    (with \( c := \codim_{P}B \geqq 0\)) over that ``base
    face''. Moreover, if \( P \) is not stout (i.e.,
    \( c > 0 \)), then the complementary simplex
    \( S_{c-1} \) is the unique minimal normally stout
    face of~\( P \).
\end{lemma}

\begin{proof}
    If the polytope is stout, then there is nothing to
    show. The general case is seen by induction on
    \( n := \dim P \geqq 2 \), with the case \( n = 2 \)
    already being settled. For \( n \geqq 3 \), we may thus
    assume that \( P \) is a pyramid \( \Pi(F) =  F * \{\mathbf{a}\} \)
    over one of its facets $F \prec P$. By induction hypothesis, the
    statement holds for that facet~\( F \). Since all faces
    of~\( P \) containing the apex~\( \ba \) are pyramids,
    every stout face already lies in~\( F \). Hence, the
    unique maximal stout face~\( B \) of~\( F \) also is the
    unique maximal stout face of~\( P \).
\end{proof}

The fundamental role played by normally stout
faces in the present approach to the HRR is that
they ``witness'' non-simplicity, cf.
\cite{brel2} 2.9:

\begin{lemma}
\label{simplemu} A polytope is simple if and
only if it has no normally stout faces.
\end{lemma}

\begin{proof} If a polytope is simple, then the links
    of all its faces are simplices, so no face is
    normally stout. On the other hand, a non-simple
    \( n \)-polytope $P$ has a vertex ${\ba} \in P$
that is incident to at least $n\,{+}\,1$ edges.
A link $L = L_P({\ba})$ of that vertex is thus
an $(n\,{-}\,1)$-polytope with more than~$n$
vertices, so it is not a simplex. Hence, as seen
above, it can be (uniquely) written as an
iterated pyramid $L = \Pi^c(B)$ (for some \( c
\geqq 0 \)) over a stout base face~$B \preceq
L$. If $c = 0$, i.e., \( L = B \), then the
vertex ${\ba}$ already is a normally stout face
of~\( P \). Otherwise we have $L = B*S_{c-1}$
with a (non-empty) simplex $S_{c-1} \precneqq L$
that is normally stout in~\( L \). Then the
unique face $F \prec P$ cutting out the face
$S_{c-1} \prec L$, i.e. such that $S_{c-1}= F
\cap L$, is normally stout in $P$.
\end{proof}

We may thus measure how ``far'' a polytope is
from being simple:

\begin{definition}
\label{defect} The defect $\mu (P) \in \NN$ of a
polytope $P$ is defined as the number of
normally stout faces of~$P$.
\end{definition}

We can restate Lemma \ref{simplemu} in these
terms: \emph{A polytope~\( P \) is simple if and
only if its defect vanishes, i.e., \( \mu(P) = 0
\).} --- Pursuing the idea sketched at the
beginning of this section, we now show that
cutting off a minimal normally stout face brings
us closer to ``simplicity'':

\begin{lemma}\label{muminus}
    Let $F \precneqq P$ be a normally stout face
    of minimal dimension, and let~\( R \) denote the
    residual polytope obtained by cutting off the
    face~\( F \) from~\( P \). Then the ``defect''
    satisfies
    \[
      \mu(R) = \mu(P) - 1 \,.
    \]
\end{lemma}

\begin{proof} No proper face $F_0 \precneqq G \cap R$ is
normally stout in $R$, since $G \cap R \prec R$,
as a cut facet, is normally trivial in $R$ and
thus $L_R(F_0) = \Pi (L_{G \cap R}(F_0))$. On
the other hand there is a bijection between the
faces of $P$ not contained in $F$ and the faces
of $R$ not contained in $G \cap R$. Since
corresponding faces have the same links and no
proper face of $F$ is normally stout in $P$ by
the minimality of $F$, we obtain $\mu (R) = \mu
(P)-1$.
\end{proof}

\begin{corollary}
\label{finite} By finitely many successive
cut-offs, every polytope is transformed into a
simple one.
\end{corollary}

\begin{proof} This follows from the above
    result together with the fact that a polytope $P$ with
    $\mu (P) = 0$ is simple, cf.
    Lemma \ref{simplemu}.
\end{proof}

So, finally, we are left with the problem to
show that the HRR for the residual polytope $R$
obtained by cutting off a minimal normally stout
face imply the HRR for the polytope $P$ itself.
To that end, we have to study the ``cut-off''
part, namely, a germ of that face. With Remark
\ref{combiprod} at our disposal, the following
result turns out to be of crucial importance,
cf. also \cite{brel2} 2.12:

\begin{lemma}
\label{normtrivmin} A normally stout face $F
\prec P$ of minimal dimension is normally
trivial and is itself a simple polytope.
\end{lemma}

\begin{proof} We let $d := \dim F$, the minimal
    dimension of any normally stout face. The case
    \( d = 0 \) being trivial, we may assume
    $d > 0$. Since an arbitrary vertex ${\ba} \in F$
    is neither simple nor normally stout in~\( P \),
its link may be written in the form $L_P({\ba})
= B * S_{c-1}$, where $B$ is stout and $c \ge
1$. The normally stout faces $F' \precneqq P$
containing ${\ba}$ correspond bijectively to the
normally stout faces of $L_P ({\ba})$ via $F'
\mapsto F' \cap L_P({\ba})$. Since $F' = F$ has
minimal dimension, and $S_{c-1}$ is the unique
normally stout face of $L_P({\ba})$ having
minimal dimension, it follows that $L_F({\ba}) =
F \cap L_P({\ba}) = S_{c-1}$, i.e., the point
${\ba}$ is a simple vertex of~\( F \), and with
\( S_{\ba} := B \precneqq L_{P}(\ba) \) in
Def.~\ref{normtrivial}, the face~$F$ is seen to
be normally trivial in~\( P \) at~${\ba}$.
\end{proof}

\section{Intersection Cohomology of Fans}
\label{intersect} In this section, we briefly
recall the construction of the intersection
cohomology of a (quasi-convex) fan $\Delta$,
referring to \cite{bbfk} or \cite{brel1} for
details. All complete fans considered in the
sequel occur as outer normal fans $\Delta(P)$
for a polytope $P \subset V$. Hence, we
systematically consider fans in the dual $V^*$
of a given vector space $V$. We are not going to
deal with non-polytopal complete fans.

\medskip
\noindent {\bf 4.A The fan space}: Motivated by
the coarse ``toric topology'' on a toric variety
given by torus-invariant open sets, we consider
a fan $\Delta$ in~$V^*$ as a finite topological
space with the subfans as open subsets. The
``affine'' fans
\begin{equation*}
  \langle \sigma \rangle \;:= \; \{\sigma\} \cup \partial\sigma
  \,\preceq\, \Delta
  \quad\hbox{with boundary fan}\quad
  \partial\sigma \;:= \; \{\tau \in \Delta \;;\; \tau \precneqq \sigma\}
\end{equation*} form a basis of the fan topology by open sets that
cannot be covered by smaller ones. Here $\preceq$ means that a
cone is a face of another cone or that a set of cones is a subfan
of some other fan. In fact, by abuse of notation, we often write
$\sigma$ instead of $\langle \sigma \rangle$, if there is no
danger of confusion.

\medskip

\noindent {\bf 4.B Sheaves}: Sheaf theory on a
fan (space) $\Delta$ is particularly simple
since a presheaf given on the basis uniquely
extends to a sheaf. In order to simplify
notation, given a sheaf $\F$ on $\Delta$, we
write
\begin{equation*}
  F_{\Lambda} :=  \F\bigl(\Lambda)
\end{equation*}
for the set of sections on the open subset (i.e., subfan)
$\Lambda \preceq \Delta$. Then a sheaf $\F$ is flabby if
and only if each restriction homomorphism
$$
\rho^\sigma_{\partial \sigma}
 \:
F_{\sigma} \to F_{\partial \sigma}
$$
is surjective.

\medskip

Here are the two most important examples:

\begin{enumerate}

\item {\bf The structure sheaf} $\A$ of~$\Delta$ is defined by
$$
A_\sigma :=   S(\affspan (\sigma)^*)\ ,
$$
the graded algebra of real-valued polynomial
functions on the subspace $\affspan (\sigma)
\subset V^*$ or rather on $\sigma$ itself, the
homomorphisms $\rho^\sigma_{\tau}: A_\sigma
\longto A_\tau$ for $\tau \preceq \sigma$ being
the restriction of functions. Hence, for
$\Lambda \preceq \Delta$, the global sections
$\in A_\Lambda$ are the $\Lambda$-conewise
polynomial functions $|\Lambda| \longto \RR$.
The grading is chosen to be twice the standard
grading, e.g. cone-wise linear functions get the
degree 2.

The structure sheaf $\A$ is flabby if and only if
$\Delta$ is a simplicial fan.

\item
  {\bf The ``equivariant'' intersection cohomology sheaf}
  $\E$ (also called ``minimal extension sheaf" in \cite{bbfk})
  is the ``smallest" flabby sheaf of graded $\A$-modules on
  $\Delta$ such that $E_\sigma$ is a finitely generated free
  $A_\sigma$-module for every cone $\sigma \in \Delta$,
  and $E_o = A_o = \RR$ for the zero cone $o :=  \{0\}$.

Let us explain the minimality condition in ``smallest'': Let
$$
  A :=  S(V) = S\bigl((V^*)^*\bigr)
$$
denote the (even-graded) algebra of polynomial
functions on the vector space $V^*$ (so in
particular, $A_\sigma = A$ for an
$n$-cone~$\sigma$, and for any fan $\Lambda$
both, $A_\Lambda$ and $E_\Lambda$ are graded
$A$-modules in a natural way). Furthermore, let
$$
  \mm :=  A^{>0}
$$
denote the unique homogeneous maximal ideal of the
graded algebra $A$. Then, given a graded $A$-module,
we define its reduction $\quer{M}$, a graded real
vector space, by
\begin{equation*}
  \quer{M} \;:= \; (A/\mm) \otimes_{A} M \, .
\end{equation*}
Since $\E$ is flabby, the reduced restriction
$$
  \qr \rho^\sigma_{\partial \sigma} \:
  \quer{E}_{\sigma} \to \quer{E}_{\partial \sigma}
$$
is also surjective. Requiring it to be even an
isomorphism means minimizing the rank of the
free $A_\sigma$-module $E_\sigma$. Note that, on
the other hand, the surjectivity of $\qr
\rho^\sigma_{\partial \sigma}$ already implies
that of $\rho^\sigma_{\partial \sigma}$.

The above conditions determine $\E$ up to isomorphy of graded
$\A$-modules, and in particular we see that $\E \cong \A$ iff
$\Delta$ is simplicial.
\end{enumerate}

\medskip

\noindent {\bf 4.C The intersection cohomology}
$IH(\Delta)$ of a complete (or more general:
"quasi-convex") fan is defined as the graded
vector space
$$
  IH(\Delta) :=   \quer{E}_\Delta\ .
$$

\medskip

\noindent {\bf 4.D Quasi-convex fans}: We call a
fan quasi-convex if it is purely
$n$-dimensional, i.e., all maximal cones are
$n$-dimensional, and the support
$|\partial\Delta|$ of its boundary subfan is a
real homology manifold or empty. Here $\partial
\Delta \preceq \Delta$ is the subfan generated
by those ($n$-1)-\allowbreak cones which are a
facet of exactly one $n$-cone in $\Delta$. In
fact, quasi-convex fans $\Delta$ are
characterized by the fact that $E_\Delta$ is a
(finitely generated) free $A$-module, cf.\
\cite{bbfk}, 4.1 and 4.4.

So in particular fans with convex or coconvex support (i.e., $V^*
\setminus |\Delta|$ is convex) as well as stars of cones in a
complete fan provide examples of such fans. Furthermore if
$\Lambda \preceq \Delta$ is a quasi-convex subfan of the complete
fan $\Delta$, we denote with $\Lambda^c \preceq \Delta$ its
(quasi-convex) complementary subfan, i.e., $\Lambda^c$ is the
subfan generated by the $n$-cones in $\Delta \setminus \Lambda$.

\medskip

\noindent {\bf 4.E Outer normal fan and
Lefschetz Operator}: Any $n$-polytope $P \subset
V$ induces a fan $\Delta = \Delta(P)$ in $V^*$
together with a strictly convex
$\Delta$-conewise linear function $\psi: V^*
\longrightarrow \RR$ as follows: For any facet
$F \preceq_1 P$ choose an "outer normal vector"
$\n_F \in V^* \setminus \{ 0 \}$, i.e., $\n_F|_F
\equiv \mathrm{ const} \ge \n_F|_P$, and denote
with $\nu (F) :=  \RR_{\ge 0}\n_F$ the
associated "outer normal ray" of the facet $F$.
To any face $G \preceq P$, we associate a cone
$\sigma (G) \subset V^*$ as follows:
$$
\sigma (G) :=  \sum_{G \preceq F \preceq_1 P} \nu (F)\ .
$$
Note in particular that $\sigma (P) = o := \{ 0
\} \subset V^*$, the zero cone. Then the {\it
outer normal fan} $\Delta (P)$ is defined as
$$
\Delta (P) :=  \{ \sigma (G); G \preceq P \}\ .
$$
We remark that $\Delta (P)$ is simplicial, iff $P$ is simple.

Denote $\bv_1,\dots,\bv_r \in V$ the vertices of
$P$. Then $\sigma_i :=   \sigma (\{\bv_i \}),
i = 1,\dots,r$ are the $n$-dimensional cones in
$\Delta$. Denote $\psi_i \in (V^*)^*$ the image
of $\bv_i$ with respect to the biduality
isomorphism $V \longrightarrow (V^*)^*$. Then
$$
\psi|_{\sigma_i} :=  \psi_i
$$
defines a strictly convex conewise linear function
$\psi \in A^2_{\Delta (P)}$. Put $\Delta := \Delta(P)$. The
multiplication map
$$
  \mu_\psi\: E_\Delta \longrightarrow E_\Delta\,,\;
  f \longmapsto \psi f
$$
induces a degree 2 map
$$
  L :=  \qr \mu_\psi \: \qr E_\Delta = IH(\Delta (P))
  \longto \qr E_\Delta  = IH(\Delta (P))\; ,
$$
the "{\bf Lefschetz operator}". We remark that
$\Delta (P + \ba) = \Delta (P)$ for $\ba \in V$
with the same Lefschetz operator, since the
correponding strictly convex functions only
differ by the ``globally linear" function $\ba
\in V \cong (V^*)^* = A^2 $.

If $\aff (P) \not = V$, the above constructions
apply \emph{mutatis mutandis} in order to give a fan
$\Delta (P)$ in $V^*/\aff_0(P)^\perp$, with the
subspace $\aff_0 (P) := \aff (P) - \ba, \ba \in
\aff (P)$, as well as a Lefschetz operator on
$IH (\Delta (P))$.

\medskip

\noindent {\bf 4.F The intersection product}:
For details cf.\ \cite{mars}. We need this
notation for a sheaf $\F$ on a quasi-convex fan
$\Delta$: The module $F_{(\Delta, \partial
\Delta)} \subset F_\Delta$ of ``sections with
compact support on $\Delta$" is defined as
$$
  F_{(\Delta, \partial \Delta)} \,:=\, \ker
  (\rho^{\Delta}_{\partial \Delta}) \, = \,
  \{ f \in F_\Delta; f|_{\partial \Delta} = 0 \}\;
  ,
$$
such that for $\Delta \preceq \Lambda$, there is
a natural inclusion $F_{(\Delta,\partial
\Delta)} \subset F_\Lambda$ by trivial extension
of sections. In order to discuss the
intersection product, we have to fix a volume
form $\omega \in \det V := \bigwedge^n V$ on
$V^*$. If the fan $\Delta$ is simplicial, we
can, following \cite{bri}, define a graded
$A$-linear ``evaluation map"
$$
  \varepsilon : A_{(\Delta,\partial \Delta)}
  \longrightarrow A[-2n]
$$
as follows: For each $n$-cone~$\sigma$, we
denote $g_\sigma \in A^{2n}_{(\sigma,
\partial \sigma)} \subset A_\sigma = A$ the unique
non-trivial
function $\ge 0$, which is the product of linear
forms in $A^2 \cong V$, whose wedge product
agrees, up to sign, with $\omega$. Then the
map~$\epsilon$ is the composite
\begin{equation}
\label{eval} E_{(\Delta,\partial \Delta)} \cong
  A_{(\Delta,\partial \Delta)} \,\subset\,
  \bigoplus_{\sigma \in \Delta^n}\! A_\sigma
  \;\longto \; Q(A) \;,\;\;
  f = (f_\sigma)_{\sigma \in \Delta^n} \longmapsto
  \sum_{\sigma \in \Delta^n}
  \frac{f_\sigma}{g_\sigma}\; ,
\end{equation}
mapping $A_{(\Delta,\partial \Delta)}$ onto $A
\subset Q(A)$. We remark that any (graded)
$A$-linear map $A_{(\Delta, \partial \Delta)}
\longrightarrow A[-2n]$ is a scalar multiple of
$\varepsilon$, and that a multiplication of
$\omega$ with a scalar $\lambda \in \RR$ results
in a multiplication of $\varepsilon$ with
$|\lambda|$.

The intersection product then is the composite
$$
.. \cap .. :  A_\Delta \times A_{(\Delta,
\partial \Delta)}
  \stackrel{\mathrm{mult}}{\longto} A_{(\Delta,
  \partial \Delta)} \stackrel{\epsilon}{\longto} A [-2n]
$$
of the multiplication of functions and the
evaluation map $\epsilon : A_{(\Delta, \partial
\Delta)} \longto A[-2n]$. In the general case
the definition uses the dual sheaf $\D\E$ of
$\E$, cf.\ \cite{mars}. Its sections over a cone
$\sigma \in \Delta$ are
$$
(\D\E)_\sigma :=  \Hom (E_{(\sigma,
\partial \sigma)},A_\sigma) \otimes \det
V_\sigma,
$$
with $V_\sigma := V/\s(\sigma)^\perp \cong \s
(\sigma)^*$. The determinant factor produces a
degree shift ($V_\sigma = A_\sigma^2$ being of
weight 2) and plays an important role in the
definition of the restriction homomorphisms
$(\D\E)_\sigma \longrightarrow (\D\E)_\tau$ for
$\tau \preceq \sigma$. Here it is necessary to
fix an orientation of $\s (\sigma)$ for every
cone $\sigma \in \Delta$, with the $n$-cones
getting the orientation defined by the volume
form $\omega \in \det V$. Then the defining
formula holds even globally:
$$
  (\D\E)_\Delta \;\cong\; \Hom \bigl(E_{(\Delta,
  \partial \Delta)},A\bigr) \otimes \det V \;\cong\; \Hom
  \bigl(E_{(\Delta, \partial \Delta)},A[-2n]\bigl)\ ,
$$
where the second isomorphy uses the isomorphism
$\det V \cong \RR, \omega  \mapsto 1$.
Furthermore there are natural isomorphisms $\E
\cong \D\E$ -- in fact, the naturality is
obtained only with the HLT for fans in lower
dimensions -- and $E_\Delta \cong
(\D\E)_\Delta$, whence we finally obtain the
intersection product
$$
.. \cap .. : E_\Delta \times E_{(\Delta,
\partial \Delta)} \longto A[-2n]\ ,
$$
which uniquely extends to a map
$$
.. \cap .. : E_\Delta \times E_\Delta \longto A
f^{-1}[-2n],
$$
where $f \in A$ is a minimal square free product
of linear forms in $A^2=V$ with $f|_{\partial
\Delta} =0$. If we apply that to the subfans
$\langle \sigma \rangle$ with an $n$-cone
$\sigma \in \Delta$ we obtain a formula
representing the intersection product of two
sections $f \in E_\Delta$ and $g \in E_{(\Delta,
\partial \Delta)}$ as a sum of local contributions:
$$
f \cap g = \sum_{\sigma \in \Delta^n} f_\sigma
\cap g_\sigma \in A
$$
with $f_\sigma:=f|_{\sigma},
g_\sigma:=g|_{\sigma}$, but note that in general
$f_\sigma \cap g_\sigma \in Q(A)$ does not
belong to $A$.

There is another way to obtain the intersection
product (cf.\ \cite{mars}, 4): Take a simplicial
refinement $\iota: \Sigma \longto \Delta$ and
realize $\E$ as a direct summand of
$\iota_*(\A)$, where $\A$ denotes the structure
sheaf of the fan $\Sigma$, (cf.\ \cite{bbfk},
2.5) -- the corresponding inclusion then is also
called a {\bf direct embedding}. Then the
composition
$$
E_\Delta \times E_{(\Delta, \partial \Delta)} \hookrightarrow
A_\Sigma \times A_{(\Sigma, \partial \Sigma)}
\stackrel{\mathrm{mult}}{\longto} A_{(\Sigma,
\partial \Sigma)} \stackrel{\epsilon}{\longto} A [-2n]
$$
of the induced embeddings and the intersection product on $\Sigma$
provides the intersection product on $\Delta$.

A third possibility is to mimic the
multiplication of functions (cf.\ \cite{mars},
4): Choose an "internal intersection product",
i.e., any symmetric $\A$-bilinear sheaf
homomorphism $\beta : \E \times \E \longto \E$
extending the multiplication of functions on the
2-skeleton -- but note that its construction
involves choices and is not natural. On the
other hand there is a distinguished section $1
\in E_\Delta$ and its image with respect to the
isomorphism
$$
E_\Delta \overset{\cong}{\longrightarrow}
(\D\E)_\Delta \cong \Hom_A (E_{(\Delta,
\partial \Delta)},A[-2n])
$$
provides an evaluation map $\varepsilon :
E_{(\Delta,\partial \Delta)} \longto A [-2n]$.
Then if we take the composite
$$
E_\Delta \times E_{(\Delta, \partial \Delta)}
\stackrel{\beta}{\longto} E_{(\Delta,
\partial \Delta)} \stackrel{\epsilon}{\longto} A
[-2n],
$$
we finally once again obtain the intersection
product!

\section{HRR for special $n$-polytopes}
\label{Perprop}
\subsection{HRR for pyramids}

\begin{proposition}
\label{HRRpyramid} If the HRR hold for polytopes
in dimension $<n$, then also for any
$n$-dimensional pyramid $P = \Pi (Q)$ over some
$(n-1)$-polytope $Q$.
\end{proposition}

\begin{proof} We may assume that $0 \in V$ is
the apex of our pyramid, i.e. $\Pi (Q)= Q *
\{0\}$. Let $\Delta:=\Delta (\Pi (Q))$ and
denote $\sigma:=\sigma (\{0\}) \in \Delta$ the
cone corresponding to the apex $0$ of the
pyramid. Then the complementary fan $\Delta_0:=
\Delta \setminus \{ \sigma \}= \langle \sigma
\rangle^c$ satisfies
$$
\Delta_0 = \mathrm{st} (\nu (Q))=\partial \sigma
+ \nu (Q):= \partial \sigma + \langle \nu (Q)
\rangle \
$$
with the outer normal ray $\nu (Q)$ of $Q
\preceq_1 P=\Pi (Q)$, and $\psi|_\sigma =0$
resp. $\psi \in A^2_{(\Delta_0, \partial
\Delta_0)} \subset A^2_\Delta$. We regard the
exact sequence
$$
0 \longto E_{(\sigma, \partial \sigma)} \longto
E_\Delta \longto E_{\Delta_0} \longto 0\ .
$$
It even splits, since $E_{\Delta_0}$ is free, the fan $\Delta_0$
being quasi-convex. Thus, there is a corresponding exact sequence
$$
0 \longto IH (\sigma, \partial \sigma) \longto IH ( \Delta)
\longto IH (\Delta_0) \longto 0\
$$
with $IH^q (\sigma, \partial \sigma) = 0$ for $q \le n$, since HLT
holds for fans in dimension $< n$, cf.\ \cite{bbfk},1.8; so for $k
\ge 0$ the restriction from $\Delta$ to $\Delta_0$ induces an
isomorphism
$$
IH^{n-k}(\Delta)
\stackrel{\cong}{\longrightarrow}
IH^{n-k}(\Delta_0) \cong IH^{(n-1)-(k-1)}(\Delta
(Q))\ .
$$
Let us comment here on the second isomorphy: The
outer normal fan $\Delta(Q)$ is a fan in
$W :=  V^*/ \RR\n_Q$, and the quotient projection
$\pi: V^* \longto W$ induces a fan map $\Delta_0
\longto \Delta (Q)$. Then, with $B :=  S(W^*)
\subset A = S((V^*)^*)$ we have
$$
E_{\Delta_0} \cong A \otimes_B E_{\Delta(Q)}\ ,
$$
whence the last isomorphism. The dual picture
looks as follows
$$
IH^{n+k}(\Delta) \cong IH^{n+k}(\Delta_0,
\partial \Delta_0) \cong IH^{n+k-2}(\Delta_0)
\cong IH^{(n-1)+(k-1)}(\Delta (Q))\ .
$$
Here the second isomorphism is the "{\bf Thom
isomorphism"}, the isomorphism induced by:
$$
E_{\Delta_0} \stackrel{\cong}{\longto}
E_{(\Delta_0.\partial \Delta_0)}, f \mapsto \psi
f\ .
$$
Analogously with $k+2$ instead of $k$ it is like
this:
$$
IH^{n+k+2}(\Delta) \cong
IH^{(n-1)+(k-1)+2}(\Delta (Q))\ .
$$
For $k > 0$, these isomorphisms transform
$$
L^k: IH^{n-k}(\Delta) \longto IH^{n+k}(\Delta)
$$
into
$$
L^{k-1}:IH^{(n-1)-(k-1)} (\Delta (Q)) \longto
IH^{(n-1)+(k-1)}(\Delta (Q))\ .
$$
This gives the HLT for $\Delta$. Now let us look at the HRR: The
homomorphism
$$
L^{k+1}: IH^{n-k}(\Delta) \longto IH^{n+k+2}(\Delta)
$$
corresponds to
$$
L^{(k-1)+1} \: IH^{(n-1)-(k-1)}(\Delta (Q)) \longto
IH^{(n-1)+(k-1)+2}(\Delta (Q))\ .
$$
So
$$
IP^{n-k}(\Delta) \cong IP^{(n-1)-(k-1)}\bigl(\Delta
(Q)\bigr) \quad \text{for}\; k
> 0\ ,
$$
while for $k = 0$ there is no contribution:
$IP^n (\Delta) =  0$ because of $L^0 = \id$. Now
the above isomorphism respects the Hodge-Riemann
forms, if we endow $V^*/\RR\n_Q$ with the volume
form $\eta$, such that $q^*(\eta) \wedge
\psi_\tau= \omega$ with the volume form $\omega$
of $V^*$, the quotient map $q: V^* \longto
V^*/\RR\n_Q$ and $\psi_\tau=\psi|_\tau \in
A^2=(V^*)^*$ with an $n$-cone $\tau \in
\Delta_0$. So the HRR hold for $\Pi (Q)$, since
they do for $Q$.
\end{proof}

\subsection{The K\"unneth formula}

We want to show that the product $S \times P_0$
of a ``HRR polytope'' \( P_{0} \) with a simple
factor~\( S \) again has the ``HRR property''.
We start with discussing the intersection
cohomology, endowed with the intersection
product.

\begin{proposition}
\label{kunneth} Let $P = S \times P_0$ be a
polytope in $V \times W$ with a simple factor~$S$, and
let $\Delta = \Sigma \oplus \Delta_0$ be the corresponding
decomposition of the respective outer normal fans. Then
there is a natural isomorphism
$$
  IH (\Delta)
  \;\stackrel{\cong}{\longrightarrow}\;
  IH (\Sigma) \otimes_{\RR} IH (\Delta_0)
$$
of graded vector spaces endowed with the intersection
forms.
\end{proposition}

\begin{proof} We let $A = S(V)$ and $B = S(W)$ denote the
    algebra of polynomials on $V^*$ and on $W^*$,
respectively. Disregarding the intersection products, the
isomorphism is seen as follows: Since~\( S \) is simple,
the fan~\( \Sigma \) is simplicial. Hence, assigning to a
cone \(\delta = \sigma \times \delta_0 \) in \( \Delta =
\Sigma \oplus \Delta_0 \) the \( A_{\delta} \)-module
$$
  E_\delta :=  A_\sigma \otimes_\RR E_{\delta_0}
$$
defines a minimal extension sheaf on $\Delta$,
as follows from an iterated application of Lemma
1.5 in \cite{bbfk}. Since the functor $A_\sigma
\otimes ...$ is exact we obtain
$$
E_{\sigma \times \Delta_0} \cong A_\sigma
\otimes E_{\Delta_0},
$$
for each $\sigma \in \Sigma$. Using the
analogous argument with the functor $... \otimes
E_{\Delta_0}$, we obtain
$$
E_\Delta \cong A_\Sigma \otimes E_{\Delta_0}.
$$
Since both, $A_\Sigma$ and $E_{\Delta_0}$, are
free modules over their base rings~\( A \) and
\( B \), respectively, the latter isomorphism
descends to the level of intersection
cohomology.

It remains to check the compatibility with the
intersection products. We first assume that the
fan $\Delta_0$ is simplicial, too. In that case,
up to suitable shifts, the tensor product of the
evaluation maps $A_\Sigma \to A$ and
$A_{\Delta_0} \to B$ defines the evaluation map
$$
  A_\Delta \,\cong\, A_\Sigma \otimes A_{\Delta_0}
  \longrightarrow A \otimes B,
$$
associated to the product of the pertinent volume
forms on $V$ and on $W$, respectively. This implies
the compatibility.

If \( \Delta_{0} \) is non-simplicial, we choose
a simplicial subdivision $\iota: \widehat
\Delta_0 \to \Delta_0$ and a direct embedding
$\E \into \iota_* (\widehat\A)$ of the
intersection cohomology sheaf $\E$ on $\Delta_0$
into the direct image of the structure sheaf
$\widehat\A$ on $\widehat \Delta_0$. It induces
a direct embedding on $\Delta = \Sigma \times
\Delta_0$.  Since these embeddings provide the
respective intersection products on
$E_{\Delta_0}$ and on $E_{\Delta}$, the
compatibility holds.
\end{proof}

To show the HRR property, we need some purely
algebraic considerations. In that framework, it
is convenient to make degrees symmetric by a
shift: Instead of \( IH(\Delta) \), graded in
even degrees ranging from \( 0 \) to \( 2n \)
and endowed with the intersection pairing and
the Lefschetz operator, we consider the
following

\begin{HRsetup}\label{HRsetup}
    Let
    \[
      W :=  \bigoplus_{k = -m}^m W^{k}
    \]
    be a finite dimensional graded vector space
    endowed with the following
    structures:
    \begin{itemize}
        \item  A non-degenerate symmetric bilinear
        form, also called the ``intersection form'',
          \[
            \left< \_\,, \_ \right> \: W \x W \longto \CC
          \]
      of total degree~\( 0 \) satisfying
      \[
        \left< W^{-k}, W^k \right> \subset i^k \RR \,,
      \]

        \item the structure of a graded
        module over the polynomial ring \( \RR[L] \)
    with \( \deg L = 2 \) such that the ``Lefschetz operator''
    $\mu_L:= L \cdot ...$ is self-adjoint
    with respect to the above form.
    \end{itemize}

    For the convenience of notation we simply
    write $L$ instead of $\mu_L$.
    These data give rise to ``HR-forms"
    $s_k (x,y) :=  (x,L^k y)$ on $W^{-k}$, and furthermore,
    to ``primitive subspaces''
    \[
      P(W^{-k}) :=  \ker(L^{k+1} \: W^{-k} \to W^{k+2}) \,.
    \]
\end{HRsetup}

\begin{definition}\label{HRmodul}
    A graded $\RR[L]$-module \( W \) endowed with
     such a structure is called an {\bf HR-module}
     if the restriction of $i^k s_k$ to the primitive
     subspace \( P(W^{-k}) \) is positive definite.
\end{definition}

We note that obviously, an HR-module satisfies $\dim W^k =
\dim W^{-k}$ (``numerical Poincar\'e duality").

The following HR-modules $A_m$ (for $m \in \NN$)
are the simple ones. They are defined by putting
$$
  A_m^k \;:= \;
  \begin{cases}
      \RR & \text{ for } -m \leqq k \leqq m
      \;\text{ with }\; k \equiv m \mod 2 \,, \\
      0 & \text{otherwise.}
  \end{cases}
 %
$$
with the intersection form mapping the pair \(
(1,1) \in A_{m}^{-k} \x A_{m}^k \) to
$\left<1,1\right> = (-i)^m$, and the ``Lefschetz
operator'' $L \: A_{m}^k \to A_{m}^{k+2}$
mapping \( 1 \mapsto 1 \), whenever that makes
sense.

\begin{remark}
\label{PropHRmodul}
\begin{enumerate}
\item[i)] Every HR-module is isomorphic to a
direct sum of modules $A_m$.

\item[ii)] A graded $\RR[L]$-submodule $U
\subset W$ of an HR-module $W$ is again an
HR-module iff it satisfies numerical Poincar\'e
duality $\dim U^{-k} =  \dim U^k$.

\end{enumerate}
\end{remark}

The link with intersection cohomology of
polytopes is provided as follows:

\begin{remark}\label{shift} Let $P$ be an
$n$-polytope and put $\Delta :=  \Delta (P)$.
Endow the graded $\RR$-vector space $W(P) :=
IH(\Delta)[-n]$ (i.e., having weight spaces
$W^k(P) :=  IH^{n+k}(\Delta)$ for \( -n \leqq k
\leqq n \)) with the intersection form,
multiplied by $(-i)^n$, and put $L$ to be the
Lefschetz operator. Then the polytope $P$
satisfies HRR iff $W(P)$ is an HR-module.
Furthermore $W(S_n) \cong A_n$ holds for the
$n$-simplex~$S_n$.
\end{remark}

We now state and prove the ``K\"{u}nneth theorem'' for
HR-modules.

\begin{proposition}
\label{HRkunneth} Let $W,W'$ be HR-modules. Then
both, $W \oplus W'$ and $W \otimes W'$ are
HR-modules, where the action of $L$ on $W
\otimes W'$ is given by $L (x \otimes y) :=   Lx
\otimes y + x \otimes Ly$.
\end{proposition}

\begin{proof} The first part of the statement
    being obvious, we only have to consider the
    tensor product. Since both, $W$ and $W'$ are
    direct sums of modules of type $A_n$ and the
    tensor product commutes with direct sums, it
    suffices to look
at $A_n \otimes A_m$\,. As a first step, a direct
computation shows that, for $n \ge 1$,
$$
  A_n \otimes A_1 \cong A_{n+1} \oplus A_{n-1}
$$
is an HR-module. By induction on $n$, it
follows that
$$
  (A_1)^{\otimes n} \cong A_n \oplus R
$$
also is an HR-module, where the ``remainder'' \(
R \) is a direct sum of terms $A_m$ with $m < n$
of the same parity as $n$. Hence the graded
vector space
$$
  B :=  A_n \otimes A_m
$$
is an $\RR[L]$-submodule of the HR-module
$C :=  (A_1)^{\otimes (n+m)}$. Satisfying numerical
Poincar\'e duality $\dim B^{-k} = \dim B^k$,
it is an HR-module itself.
\end{proof}

\begin{corollary}
\label{kunnethconclud} If, in the situation of
Proposition \ref{kunneth}, the polytope $P_0$
satisfies HRR, then so does $P = S \times P_0$.
\end{corollary}

\begin{proof} According to
\ref{kunneth} and in the notation of
\ref{shift}, the graded $\RR [L]$-module $W(P)$
can be written as
$$
  W(P) \cong W(S) \otimes W(P_0)\,;
$$
hence the claim follows from \ref{HRkunneth} and
the HRR for the simple polytope $S$.
\end{proof}

\subsection{Transversal Cuttings}

\begin{proposition}
\label{cutoff} If the affine hyperplane $H$ cuts
$P$ transversally into the polytopes $P_1$ and
$P_2$, i.e., $H$ has nonempty intersection with
the relative interior of $P$ and does not contain
vertices of $P$, then the validity of HRR for $P_1$
and $P_2$ and for lower dimensional polytopes
implies HRR for $P$.
\end{proposition}

\begin{proof}
Let us write
$$
F_i :=   P \cap H \preceq P_i,
$$
using the index $i = 1,2$ in order to indicate
when $P \cap H$ should be considered as a facet
of $P_i$.

\medskip

\noindent {\bf (A) Fans involved.} First of all
let, as usual,
$$
  \Delta :=  \Delta (P) \;, \;\; \text{and put}\;\;
  \Delta_i :=  \Delta(P_i) \;\;\text{for}\;\; i = 1,2\,.
$$
Secondly, we consider the ``intermediate'' polytope
$Q$ cut out from $P$ by $H$ and a nearby parallel hyperplane.
Its outer normal fan \( \Delta(Q) \) is obtained by just
putting together the stars
$$
  \Lambda_i :=
  \mathrm{st}(\nu (F_i)) \preceq \Delta_i
$$
of the outer normal ray to the ``cut'' facet with respect to
the fans \( \Delta_{i} \) to  a new complete fan
$$
  \Lambda :=  \Lambda_1 \cup \Lambda_2 = \Delta (Q) \,.
$$
Finally, we let $\Phi :=  \Delta (F_1) = \Delta (F_2)$ denote
the outer normal fan of $F :=  F_1 = F_2$ in $W :=
V^*/(\RR \cdot \n_{F_i})$ (note that $\n_{F_2} = -\n_{F_1}$).

\medskip

\noindent {\bf (B) Gluing of $\boldsymbol{IH}$.} Let $\G$
and $\E$ denote the respective intersection cohomology
sheaves on $\Phi$ and $\Lambda$.  The projection $\pi :
V^* \to W$ induces a map of fans $\Lambda \to \Phi$, and
we have $\E \cong \pi^*(\G)$.
In particular, writing $B :=   S(W^*) \subset A :=   S((V^*)^*)$,
there is a natural injection
$$
  A \otimes_B G_{\Phi} \hookrightarrow E_{\Lambda},
$$
which after restriction to the subfans $\Lambda_i \preceq \Lambda$
gives isomorphisms
$$
  E_{\Lambda_1} \stackrel{\cong}{\longleftarrow} A \otimes_B
  G_{\Phi} \stackrel{\cong}{\longrightarrow} E_{\Lambda_2}.
$$
Denote with
$$
  S: E_{\Lambda_1} \stackrel{\cong}{\longto} E_{\Lambda_2}
$$
the resulting $A$-module isomorphism. Now consider the exact
sequence
\begin{equation}
\label{Gl1} 0 \longrightarrow K \longrightarrow E_{\Delta_1}
\oplus E_{\Delta_2} \longrightarrow E_{\Lambda_2} \longrightarrow
0,
\end{equation}
where the second nontrivial map sends
$(f_1,f_2)$ to $S(f_1|_{
\Lambda_1})-f_2|_{\Lambda_2}$, and $K \subset
E_{\Delta_1} \oplus E_{\Delta_2}$ is the kernel
of that map. Since $E_{\Lambda_2}$ is a free
$A$-module and
$$
IH(\Lambda_i) \cong IH(\Phi),
$$
the exact sequence (\ref{Gl1}) induces the exact sequence
\begin{equation}
\label{Gl1R} 0 \longrightarrow \overline K \longrightarrow IH
(\Delta_1) \oplus IH(\Delta_2) \longrightarrow IH (\Phi)
\longrightarrow 0.
\end{equation}
Furthermore we need
\begin{equation}
\label{Gl2} 0 \longrightarrow D \longrightarrow K \longrightarrow
E_\Delta \longrightarrow 0,
\end{equation}
where the second nontrivial map is gluing of
sections: The fan $\Delta$ is the union
$$
\Delta =  \Lambda^c_1 \cup \Lambda^c_2
$$
of the complementary subfans $\Lambda^c_i \preceq \Delta_i$ of
$\Lambda_i \preceq \Delta_i$, i.e.,
$$
\Delta_i =  \Lambda_i \cup \Lambda^c_i
$$
with the quasi-convex fans $\Lambda_i$ and $\Lambda^c_i$
intersecting only in their common boundary fan and, in particular,
$|\Lambda_i| = |\Lambda^c_j|$ for $j \not =  i$. Now for a pair
$(f_1.f_2) \in K \subset E_{\Delta_1} \oplus E_{\Delta_2}$ we
define its image as the section $f \in E_\Delta$ satisfying
$f|_{\Lambda^c_i} :=  f_i|_{\Lambda^c_i}$. The kernel $D$ then
satisfies
$$
D = \{ (h,S(h)); \ h \in E_{(\Lambda_1,\partial \Lambda_1)} \} \cong
E_{(\Lambda_1,\partial \Lambda_1)}\ .
$$

\bigskip

\noindent {\bf (C) Gluing of the intersection
product.} The exact sequence (\ref{Gl2}) yields an
isomorphism
$$
E_\Delta \cong K/D,
$$
in fact that quotient representation holds even
with respect to the intersection pairings on
$E_\Delta$ and $K \subset E_{\Delta_1} \oplus
E_{\Delta_2}$: Consider two pairs $(f_1,f_2),
(g_1,g_2) \in K$ and denote $f,g \in E_\Delta$
their respective images in $E_\Delta$. Then in
$Q(A)$ we obtain
$$
(f_1,f_2) \cap (g_1, g_2) = f_1 \cap g_1 + f_2
\cap g_2
$$
\begin{equation}
\label{IS}  = f_1 \cap_{\Lambda_1^c} g_1 +
f_1\cap_{\Lambda_1} g_1 + f_2 \cap_{\Lambda_2}
g_2 + f_2 \cap_{\Lambda_2^c} g_2
\end{equation}
$$
 = f_1 \cap_{\Lambda_1^c} g_1 + f_2
\cap_{\Lambda_2^c} g_2 = f \cap g,
$$
since the middle terms in (\ref{IS}) add up to 0.
This can be seen as follows: The restrictions
$E_{\Delta_i} \longrightarrow E_{\Lambda_i}$,
$i = 1,2$, combine to a map
$$
E_{\Delta_1} \oplus E_{\Delta_2} \supset K \longrightarrow
E_\Lambda \subset E_{\Lambda_1} \oplus E_{\Lambda_2}
$$
with image $A \otimes_B G_{\Phi}$. Denote $\hat
f, \hat g \in E_\Lambda$ the respective images
of the pairs $(f_1,f_2), (g_1,g_2) \in K$. Then
$$
\hat f \cap \hat g = f_1\cap_{\Lambda_1} g_1 +
f_2 \cap_{\Lambda_2} g_2.
$$
So what we finally have to prove is that
$$
A \otimes_B G_{\Phi} \subset E_{\Lambda}
$$
is an isotropic subspace. For this we may even
assume that $\Phi$ and thus also $\Lambda$ is
simplicial (remember that $\Phi \cong \partial
\Lambda_i$), otherwise take a simplicial
refinement $\iota: \widehat \Phi \longto \Phi$
and a direct embedding $\G \longto \iota_*(\A)$,
where $\A$ is the structure sheaf of $\widehat
\Phi$. There is an induced simplicial refinement
$\widehat \Lambda \longrightarrow \Lambda$ and
embedding $E_\Lambda \hookrightarrow A_{\widehat
\Lambda}$ respecting the intersection product.
Since $A \otimes_B G_\Phi \subset A \otimes_B
A_{\widehat \Phi}$, our claim holds for
$\Lambda$, if it does for $\widehat \Lambda$.

So let us now consider the case where $\Phi$ and
hence also $\Lambda$ is simplicial. In that
situation, it suffices to check that the
evaluation map $\epsilon : E_{\Lambda} =
A_{\Lambda} \longto A[-2n]$ vanishes on $A
\otimes_B G_{\Phi}$. We use the formula
(\ref{eval}) in section 4.F for $\varepsilon$:
$$
  E_{\Lambda} \cong
  A_{\Lambda} \,\subset\,
  \bigoplus_{\sigma \in \Lambda^n}\! A_\sigma
  \;\longto \; Q(A) \;,\;\;
  f = (f_\sigma)_{\sigma \in \Lambda^n} \longmapsto
  \sum_{\sigma \in \Lambda^n}
  \frac{f_\sigma}{g_\sigma}\; .
$$
But the $n$-cones in $\Lambda$ may be grouped in
pairs $\sigma_i \in \Lambda_i, i = 1,2$ with $\pi
(\sigma_1) = \pi (\sigma_2)$. Then for $f \in A
\otimes_B G_{\Phi} = A \otimes_B A_{\Phi}$, we have
$$
  f_{\sigma_1} = f_{\sigma_2}\ , \ \mathrm{while}\ \
  g_{\sigma_1} = -g_{\sigma_2}\ ,
$$
and thus $\epsilon (f) = 0$. ---

\medskip
\noindent {\bf (D) The HR relations for $\boldsymbol{P}$.}
Since
$$
\overline D \cong IH^*(\Lambda_i, \partial \Lambda_i) \cong
IH^{*-2}(\Lambda_i) \cong IH^{*-2}(\Phi)
$$
and the third module $E_\Delta$ in the exact
sequence (\ref{Gl2}) is free, there is an
associated exact sequence
$$ 0 \longrightarrow IH^{*-2}(\Phi) \longrightarrow \overline K
\longrightarrow IH (\Delta) \longrightarrow 0,
$$
realizing $IH (\Delta) =  \overline E_\Delta$ in the same way as
before $E_\Delta$. Furthermore it is compatible with the natural
Lefschetz operators on all three terms; we shall denote them
simply $L$ in all cases.

Now let $\zeta \in IH^{n-k} (\Delta)$, $\zeta \not =
0$, be a primitive class, i.e., $L^{k+1}(\zeta) = 0$. We
can lift it to a pair $\xi = (\xi_1,\xi_2) \in
\overline K$. We may actually assume that the classes
$\xi_i \in IH^{n-k}(\Delta_i)$ are again primitive, i.e.,
$\xi_i \in IP^{n-k}(\Delta_i)$: Because of
$L^{k+1}(\zeta) = 0$, we have
$$
  L^{k+1}(\xi) = \eta \in IH^{n+k}(\Phi)
  \subset \overline K^{n+k+2} \,,
$$
and since
$$
L^{k+1}: IH^{n-k-2}(\Phi) =
IH^{(n-1)-(k+1)}(\Phi) \stackrel{\cong}{\longrightarrow}
IH^{(n-1)+(k+1)}(\Phi) = IH^{n+k}(\Phi)
$$
is an isomorphism -- by the assumption, HRR and
thus HLT holds for the lower dimensional
polytope $P_1 \cap P_2$
-- , we may replace $\xi$ with $\xi-
L^{-(k+1)}(\eta)$.

So now let both $\xi_1$ and $\xi_2$ be
primitive. Since $\zeta \not=  0$, we have
$\xi_i \not= 0$ for at least one index $i$. Thus
$$
(-1)^{(n-k)/2} \zeta \cap L^k
(\zeta) = (-1)^{(n-k)/2} \xi_1 \cap L^k
(\xi_1) + (-1)^{(n-k)/2} \xi_2 \cap L^k
(\xi_2)> 0,
$$
since the HRR hold for $P_1$ and $P_2$.
\end{proof}

\section{Deformation}
\label{Deformation}

\begin{proposition}
\label{deform} The germ $G = G_P(F)$ of a
normally trivial face $F \prec P$ can be
deformed into the product $F \times \Pi (L)$,
where $L :=  L_P(F)$ denotes a link of $F$ in
$P$.
\end{proposition}

\begin{proof}
We use the terminology of \ref{combiprod} and
assume that $\overset{\circ}{F} \cap N =  \{ 0
\}$; so the affine span $U :=  \aff (F)$ as well as
$N$ are linear subspaces, and $V =  U \oplus N$.
Furthermore write $N =  W \oplus \RR$, such that
$H = U \times W \times \{1\}$, and $L_P(F) = L
\times \{1\}$ with a polytope $L \subset W$.
Denote $\bu_1,\dots,\bu_r \in U$ the vertices of
$F$ and $\bw_1,\dots,\bw_s \in W$ the vertices of
$L$, then $G$ has vertices $(\bu_i,0,0)$ and
$(\bu_i+\bu_{ij},\bw_j,1)$ with suitable vectors
$\bu_{ij} \in U$.

Now let $G_t :=  G \cap (U \times W \times [0,t])$ be the truncated
germ, and consider on
$$ V  =  U \oplus N  =  U
\oplus (W \oplus \RR)
$$
the linear isomorphism
$$
F_t :=   \id_U \oplus t^{-1} \id_N.
$$
Then $(0,1] \ni t \mapsto Q_t :=  F_t (G_t)$
extends to a deformation $[0,1] \ni t \mapsto Q_t$
with $Q_0 =  F \times \Pi (L)$. In fact the polytope
$Q_t$ has vertices $(\bu_i,0,0)$ and
$(\bu_i+t\bu_{ij},\bw_j,1)$.
\end{proof}

\begin{theorem}
\label{HRRgerm} Let $F \prec P$ be a normally
trivial face of an $n$-dimensional polytope $P$,
and assume that $F$ itself is a simple polytope.
Then the HRR hold for $G = G_P(F)$ if they hold
for lower dimensional polytopes.
\end{theorem}

\begin{proof}
For $\dim F =  0$, the germ $G = \Pi (L_P(F))$
is a pyramid, so we may apply Proposition
\ref{HRRpyramid}. Now let $\dim F > 0$. Let us
first give a

\bigskip

\noindent {\bf Survey of proof}: We consider the
deformation $Q_t, 0 \le t \le 1$, of Prop.\
\ref{deform}, which deforms $Q_1 = G$ into $Q_0
=  F \times \Pi (L)$. The HLT and the HRR hold
on $Q_0$ according to the K\"unneth formula.
Then we show that the HLT holds on $Q_t$ for all
$t \in (0,1]$, cf.\ \ref{hltdeform}, hence the
HR-forms on $IH (\Delta (Q_t))$ are
non-degenerate for any $t \in [0,1]$. Since the
combinatorial type of the polytopes $Q_t$ is
constant along the deformation, the Betti
numbers are so too, in fact both $IH^{n-k}
(\Delta (Q_t)), t \in [0,1]$ and the $k$-th
HR-form $s_k^t$ on $IH^{n-k} (\Delta (Q_t))$
depend continuously on $t \in [0,1]$, cf.\
\ref{continuity}. Since they are non-degenerate,
their signature is independent of the parameter
$t \in [0,1]$. Hence the HR-equations
\ref{HRequat} hold for all $t$, since they do
for $t = 0$.

\bigskip

\noindent {\bf (A) The deformation on the fan level}: \noindent

\noindent {\bf The case $t > 0$}: For $t >0$ the fan $\Delta_t :=
\Delta (Q_t)$ is ``linearly equivalent'' to $\Delta_1 =  \Delta
(G)$, i.e., $\Delta_t$ is the image of $\Delta_1$ with respect to a
linear isomorphism of the vector space $V^*$, namely the inverse
$(F_t^*)^{-1}$ of the dual $F_t^*$ of the map $F_t : V
\longrightarrow V$ transforming $G_t$ into $Q_t$ in the proof of
\ref{deform}. It provides an isomorphism
$$
F_t^*: \Delta_t = \Delta (Q_t) \longrightarrow \Delta (G_t) = \Delta
(G) = \Delta (Q_1)  =  \Delta_1\ .
$$
{\bf Behaviour near $0$}: We replace the linear isomorphism
$(F_t^*)^{-1}$ mapping $\Delta_1$ onto $\Delta_t$ with a
$\Delta_0$-conewise-linear isomorphism $S_t : V^* \longrightarrow
V^*$, such that
$$
S_t (\Delta_0) =  \Delta_t\ .
$$
The construction of $S_t$ is as follows: Consider the subfan
$\Gamma \preceq \Delta (G)$ generated by the cones $\sigma (F_0)$,
where $F_0 \preceq G$ is a minimal face projecting onto the entire
pyramid $\Pi(L)$, i.e., $\pi (F_0) = \Pi (L)$ with the projection
$\pi : V =  U \oplus N \longrightarrow N$. The support $|\Gamma|$ is
the graph of a map $H:U^* \longrightarrow N^*$. In fact that map
is $\Phi$-conewise linear for the outer normal fan $\Phi :=  \Delta
(F)$ of the polytope $F \subset U$, and
$$
S_t: U^* \oplus N^* \longrightarrow U^* \oplus N^*, \ (\bx,\by)
\mapsto (\bx ,\by+ t H(\bx))
$$
then defines a $\Delta_0$-conewise linear isomorphism with the
desired properties. Note here that $\Delta_0$ is the product
$$
\Delta_0 =  \Phi \times \Lambda
$$
of the (simplicial) fan $\Phi :=  \Delta (F)$ in $U^*$ and the fan
$\Lambda :=   \Delta(\Pi (L))$ in $N^*$.

\medskip

\noindent{\bf (B) Pull back isomorphisms}: Both $(F_t^*)^{-1}$ and
$S_t$ act on the global sections of the structure resp. the
intersection cohomology sheaf by pull back. Let us write
$$
A_t :=  A_{\Delta_t}\ ,\ E_t :=  E_{\Delta_t}.
$$
Then $(F_t^*)^{-1}$ induces $A$-module
isomorphisms
$$
A_t \longrightarrow A_1\ ,\ E_t \longrightarrow E_1\ ,
$$
in particular $IH(\Delta_t) \cong IH(\Delta_1)$ in a natural way,
while for $S_t$ the corresponding maps
$$
A_t \longrightarrow A_0\ ,\ E_t \longrightarrow E_0\ ,
$$
both denoted $S_t^*$, are only isomorphisms of graded vector
spaces due to the fact that for the subalgebra $A \subset A_t$ of
"global polynomials" we in general have $S_t^*(A) \not\subset A
\subset A_0$. So we can not any longer identify $IH(\Delta_t)$ in
a reasonable way with $IH(\Delta_0)$.
\bigskip

\noindent {\bf (C) Continuity statements}: That
everything is continuous in $t \in (0,1]$
follows now immediately from the fact that the
strictly convex function on $\Delta_1 =  \Delta
(G)$ given by the vertices of $G_t$ is
continuous in $t$. Near $0$ there is no natural
trivialization of the family $IH (\Delta_t)$;
instead we have to represent $IH(\Delta_t)$ as a
factor space $E/M_t$ of a bigger vector space
$E$ independent of $t$ with varying subspace
$M_t$:

\begin{proposition}
\label{continuity} There is a finite dimensional graded vector
space $E$ and continuous families of
\begin{enumerate}
\item subspaces $M_t \subset E$ of constant
dimension, such that in a natural way
$$
IH (\Delta_t) \cong E/M_t\ ,
$$
\item endomorphisms $\hat L_t:E \longrightarrow
E$ with $\hat L_t (M_t) \subset M_t$ inducing
the Lefschetz operator of $\Delta_t = \Delta
(Q_t)$,

\item symmetric bilinear forms $\beta_t :E \times E
\longrightarrow \RR$ with $\beta_t (M_t,E) = 0$, inducing the
intersection product.
\end{enumerate}
\end{proposition}

\begin{proof} Let us start with
\medskip

\noindent {\bf The vector spaces $E$ and $M_t$}: Write
$\Delta :=  \Delta_0$. We take:
$$
E :=  \bigoplus_{q = 0}^{2n} E^q_0
$$ and
$$
M_t :=   S_t^* (\mm E_t) \cap E \subset E.
$$
The subspaces $M_t$ can be represented in the
form $\Phi_t (\mm^{< 2n} \otimes E)$ with the
continuous family of linear maps
$$
\Phi_t : \mm^{< 2n} \otimes E \longto E,\ g
\otimes f \mapsto S_t^*(g)f;
$$
furthermore, since $\Delta_t$ and $\Delta$ are
combinatorially equivalent, we get that $\dim
M_t =  \dim E - \dim IH (\Delta_t)$ is
independent of $t$. The map $\hat L_t$ is
multiplication with $S_t^* (\psi_t) :=   \psi_t
\circ S_t$, except on the highest weight
subspace $E_\Delta^{2n}$, where it vanishes.
Here $\psi_t$ denotes the strictly convex
function belonging to $Q_t$.

\medskip

\noindent {\bf Continuity of the intersection
product}: Here we consider in general the
situation, where we have a fan $\Delta$ in $V^*$
and $\Delta_t :=  S_t (\Delta)$ with a continuous
family of $\Delta$-conewise linear isomorphisms
$S_t: V^* \longrightarrow V^*$. The bilinear
form we consider is
$$
\beta_t: E \times E \subset E_\Delta \times E_\Delta
\stackrel{(S_t^{-1})^*}{\longto} E_t \times E_t
\stackrel{\cap}{\longrightarrow} A[-2n] \longrightarrow A^0 \cong
\RR
$$
where the last arrow is (up to the shift) the projection $A =  A^0
\oplus \mm \longrightarrow A^0$.

 Let us first look at the case of a

\medskip

\noindent {\bf Simplicial fan $\Delta$}: Then we
have $\E = \A$. Take $r,s$ with $r+s = 2n$. We have
to show that the map
$$
A_0^r \times A_0^s
\stackrel{(S_t^{-1})^*}{\longto} A^r_t \times
A^s_t \stackrel {\cap}{\longto} A^0  =  \RR
$$
depends continuously on $t \in [0,1]$. But that map may be
rewritten as
$$
A^r_0 \times A^s_0
\stackrel{\mathrm{mult}}{\longto} A_0^{2n}
\stackrel{(S_t^{-1})^*}{\longto} A^{2n}_t
\stackrel{\epsilon_t}{\longto}   A^0  = \RR \ ,
$$
using the fact that $(S_t^{-1})^*$ commutes with
the multiplication of functions. Here
$\epsilon_t$ denotes the restriction of the
evaluation map
$$
A_t \longto A[-2n]\
$$
to the $2n$-th weight space $A_t^{2n}$. It is
well defined after having fixed a volume form on
$V$. So, eventually we have to check that the
map
$$
A_0^{2n} \stackrel{(S_t^{-1})^*}{\longto} A^{2n}_t
\stackrel{\epsilon_t}{\longto}  \RR
$$
depends continuously on $t \in [0,1]$. Take any
$n$-cone $\sigma_0 \in \Delta  =  \Delta_0$, set
$\sigma_t :=  S_t (\sigma_0) \in \Delta_t$ and
choose a non-negative function $f_t \in
A^{2n}_{(\sigma_t,
\partial \sigma_t)} \subset A^{2n}_t$, the product of
linear forms $\in (V^*)^*$, whose
$\wedge$-product is up to sign the volume form
on $V$. Denote $T_t: V^* \longrightarrow V^*$
the linear map, which coincides with $S_t$ on
$\sigma_0$. Then we have $S_t^*(f_t) = \det
(T_t)\, f_0$, and the map $\epsilon_t \circ
(S_t^{-1})^*: A^{2n}_0 \longto \RR$ can be
thought of as $\det (T_t)^{-1}$ times the
projection operator $A^{2n}_0 \longto \RR f_0
\subset A^{2n}_0$ with kernel $M_t^{2n} \subset
E^{2n} = A_0^{2n}$, since $\epsilon_t (f_t) =
1$. That yields the desired continuity with
respect to $t \in [0,1]$.

\medskip

\noindent {\bf The general case}: Take a simplicial refinement
$\Sigma \stackrel{\iota}{\longrightarrow} \Delta$ and consider an
enbedding
$$
\E \hookrightarrow \iota_*( \A) \ \ \mathrm{and}\ E_\Delta
\hookrightarrow A_\Sigma
$$
as in section 4.F. For $\Sigma_t :=  S_t (\Sigma)$ it induces
embeddings
$$
\label{emb} E_t :=   E_{\Delta_t} \hookrightarrow A_t :=  A_{\Sigma_t}\
,
$$
which according to \cite{mars} respect the intersection pairings.
Then the intersection product takes the form
$$
E_0^r \times E_0^s \longrightarrow A_0^r \times
A_0^s \stackrel{(S_t^{-1})^*}{\longto} A^r_t
\times A^s_t \stackrel {\cap}{\longto} A^0  =
\RR\ ,
$$
so the simplicial case applies.
\end{proof}
\bigskip

\noindent {\bf (D) Statement HLT for the $Q_t, t
\in (0,1]$}: Finally we show
\begin{proposition}
\label{hltdeform} The HLT holds for the
polytopes $Q_t, t \in (0,1]$, in particular the
HR-forms are non-degenerate for all $t \in
[0,1]$.
\end{proposition}

\begin{proof} Since the linear isomorphism $F_t^*: V^* \longrightarrow
V^*$ induces an isomorphism
$$
IH (\Delta_t) \stackrel{\cong}{\longrightarrow} IH(\Delta_1)
$$
we may replace $\Delta_t$ with
$\Delta :=  \Delta_1 =  \Delta (G) = \Delta (G_t)$.
Denote $\psi :=  \psi_t$ the function given by the
vertices of $G_t$ and $L = L_t$ the corresponding
Lefschetz operator. Because of Poincar\'e
duality it suffices to prove that the
corresponding $k$-th iteration of the Lefschetz
map $L^k : IH^{n-k}(\Delta) \longto
IH^{n+k}(\Delta)$ is injective for $k>0$. We
have
$$
\Delta = \Theta \cup \Theta_0\ ,
$$
with the subfan $\Theta =
\mathrm{st}(\sigma(F))$ corresponding to the
ridge $F \preceq G$ of the "hip roof" $G$, and
the subfan $\Theta_0 :=  \mathrm{st}(\nu (G \cap
H))$ associated to the bottom or cut facet $G
\cap H$. The rays of $\Delta$ not contained in
$\sigma (F)$ are the outer normal rays $\rho_i
:=  \nu (F_i)$ of the facets $F_0,\dots,F_r$ of
$G$ not containing $F$ - say, $F_0 :=  G \cap H
\prec G$ is the bottom of $G$. Since $F$ is
simple, they are normally trivial in $G$ (To see
that use the fact that $G$ is combinatorially
equivalent to $F \times \Pi (L)$, cf.\
\ref{combiprod}, and that the $F_1,\dots,F_r$
under that equivalence correspond to facets of
$F$ times the pyramid $\Pi (L)$, while for the
cut facet $F_0$ the claim is obvious), in
particular we find that any $n$-cone $\sigma
\succeq \rho_i$ is the sum of $\rho_i  = \nu
(F_i)$ and the unique opposite facet of $\sigma$
(corresponding to the unique edge starting in
the vertex the cone $\sigma$ is associated with,
and not contained in the facet $F_i$.)

As a consequence, there are (unique) functions $\psi_i \in
A^2_\Delta$ vanishing outside $\Theta_i :=  \mathrm{st} (\nu (F_i))$
with $\psi_i =  \psi$ on the ray $\nu (F_i)$ for $i = 0,\dots,r$. On the
other hand, we may assume $0 \in F$ resp. $\psi|_{\sigma (F)} = 0$.
So altogether we have
$$
\psi  =  \sum_{i = 0}^r \psi_i\ .
$$
Now assume $k \equiv n \ $mod$(2), 0 < k \le n$ and $\xi \in
IH^{n-k}(\Delta)$ with $L^k(\xi) = 0$. For
$\Theta_i :=  \mathrm{st}(\nu (F_i))$ we show $0 =  \xi|_{\Theta_i} \in
IH (\Theta_i)$. First of all $L^k(\xi) = \psi^k \xi = 0$ gives
$$
0 =  \xi \cap \psi^k\xi   =  \sum_{i = 0}^r \xi \cap \psi_i
\psi^{k-1}\xi  =  \sum_{i = 0}^r  \xi_i \cap \phi_i ^{k-1}\xi_i
$$
with $\xi_i :=  \xi|_{\Theta_i} \in
IH^{n-k}(\Theta_i) \cong IH^{(n-1)-(k-1)}(\Delta
(F_i))$, the strictly convex conewise linear
function $\phi_i \in A^2_{\Delta (F_i)}$ being
the pullback of $\psi|_{\partial \Theta_i}$ with
respect to the inverse of $\pi|_{\partial
\Theta_i}: \partial \Theta_i \stackrel {\cong}
{\longto} \Delta (F_i)$, where $\pi$ is the
quotient projection $\pi:V^* \longrightarrow
V^*/\s (\nu (F_i))$, the intersection product
$\xi_i \cap \phi_i^{k-1}\xi_i$ referring to
$\Delta (F_i)$. The equality $\xi \cap \psi_i
\psi^{k-1}\xi  =  \xi_i \cap \phi_i ^{k-1}\xi_i$
is obtained as follows: If we use the internal
product approach to the intersection product, we
obtain a commutative diagram
$$
\begin{array}{ccc}
    (\xi_i, \xi_i)     & & (\xi, \psi_i \xi)|_{\Theta_i}  \\
 \in & & \in  \\
    IH ^{(n-1)-(k-1)}(\Delta(F_i)) \times IH^{(n-1)-(k-1)}
(\Delta(F_i)) & \longrightarrow & IH^{n-k} (\Theta_i)
\times IH^{n-k+2} (\Theta_i, \partial \Theta_i) \\
\downarrow & & \downarrow \\
IH^{2n-2}(\Delta(F_i)) & \stackrel{\cong}{\longrightarrow} &
IH^{2n} (\Theta_i,
\partial \Theta_i)
\end{array},
$$
where the left vertical arrow denotes the
($k$--1)rst HR-form of $\Delta (F_i)$ and the
right one is the intersection product composed
in one argument with multiplication with
$\psi^{k-1}$. The second component of the upper
horizontal homomorphism and the lower one are
pull back to the fan $\Theta_i$ followed by the
Thom isomorphism (cf.\ the proof of
\ref{HRRpyramid}), i.e., multiplication with
$\psi_i$.

Now as a consequence of $L^k (\xi) = 0$ we have
$\xi_i \in IP^{(n-1)-(k-1)} (\Delta (F_i))$.
Hence the HRR for $\Delta (F_i)$ give that all
the summands are either non-negative or
non-positive. So necessarily $\xi_i  \cap
\phi_i^{k -1}\xi_i  = 0$ resp. $\xi_i = 0$ for
$i = 0,\dots,r$. In particular $\xi_0 = 0$ and
 the exact sequence
$$
0 \longrightarrow IH^{n-k}(\Theta, \partial \Theta)
\longrightarrow IH^{n-k}(\Delta) \longrightarrow
IH^{n-k}(\Theta_0) \longrightarrow 0\
$$
tells us that $\xi \in IH^{n-k}(\Theta, \partial \Theta) \subset
IH^{n-k}(\Delta)$, and thus, in order to conclude $\xi = 0$, it
suffices to prove that $L^k|_{IH^{n-k}(\Theta, \partial \Theta)}$
is injective for $k
> 0$ resp. that dually $L^k: IH^{n-k}(\Theta) \longrightarrow
IH^{n+k}(\Theta)$ is surjective for $k > 0$. In fact we show, that
the graded vector space $IH (\Theta) /L^k(IH (\Theta))$ has
weights at most $n+k-1$.

Now the fan $\Theta$ has the form
$$
\Theta =  S (\Phi \times  \sigma (0))
 \subset U^*
\oplus N^*
$$
with the cone $\sigma (0) \subset N^*$
associated to the apex $0 \in \Pi (L)$ of the
pyramid $\Pi (L) \subset N$ and $S :=  S_1$. Using
the induced vector space isomorphism
$$
S^* : E_{\Theta} \longto E_{\Phi \times  \sigma
(0)} \cong A_\Phi \otimes E_{\sigma (0)}
$$
we can write
$$
IH(\Theta) \cong A_\Phi \otimes E_{\sigma (0)}/
S^*(\mm)(A_ {\Phi} \otimes E_{\sigma (0)}))\
$$
with $S^*(\mm) \subset A_{\Phi \times  \sigma
(0)}$. Since $\psi|_{\sigma (F)}  =  0$ (because
of $0 \in F$) and $\Theta =  \mathrm{st} (\sigma
(F))$ (remember that $\sigma (F) =  o \times
\sigma (0) \subset U^* \oplus N^*$), we can
write $\psi|_{\Theta} =  \chi \circ p$ with the
projection $p:U^* \oplus N^* \longto U^*$ and a
function $\chi \in A^2_\Phi$. Now in order to
compute $IH (\Theta)/L^k(IH (\Theta))$ we have
to regard on $A_\Phi \otimes E_{\sigma (0)}$ the
"twisted" $A$-module structure obtained from
that of $E_{\Theta}$ by pull back via $S$, with
other words, a function $f \in A$ acts on
$A_\Phi \otimes E_{\sigma(0)}$ by "standard"
multiplication with $f \circ S \in A_{\Phi
\times \sigma (0)}$.

Now write $A = C \otimes D$ with the polynomial
algebras $C :=   S((U^*)^*)$ resp. $D :=   S((N^*)^*)$
on $U^*$ resp. $N^*$. Then $g = g (\bx) \in C$
acts on the first factor $A_\Phi$ only, while $h
\in D$ acts by standard multiplication with
$h(\by + H(\bx))$. We have now to divide by the
submodule obtained by multiplication with
$\mm_C, \chi^k, \mm_D$. Looking first at $\mm_C$
and $\chi^k$ gives
$$
(IH(\Phi)/L^k(IH(\Phi)) \otimes E_{\sigma(0)}\ ,
$$
a $D$-module. The graded vector space $F :=
IH(\Phi)/L^k(IH(\Phi))$ has weights at most $s
+k-1, \ s :=  \dim U^* < n$, according to the
HLT for $F$. It admits a descending filtration
by the $D$-submodules
$$
F^{\ge i} \otimes E_{\sigma(0)}, \ 0 \le i \le
s+k\ ,
$$
with free successive quotients
$$
  F^{\ge i} \otimes E_{\sigma(0)}/ F^{\ge i+1}
  \otimes E_{\sigma(0)} \;\cong\; (F^{\ge i}/F^{\ge
i+1}) \otimes E_{\sigma (0)}\ ,
$$
since $h \in D$ acts only on the second factor
of the right hand side --- the twist being
factored out. The short exact sequences
$$
0 \longrightarrow F^{\ge i+1} \otimes E_{\sigma
(0)} \longrightarrow F^{\ge i} \otimes E_{\sigma
(0)}\longrightarrow F^{\ge i} \otimes E_{\sigma
(0)}/ F^{\ge i+1} \otimes E_{\sigma (0)}
\longrightarrow 0\
$$
remain exact after reduction mod $\mm_D$: The
third terms being free $D$-modules, they are
split. So, since the third non-trivial term has
weights $< i+t$ with $t := \dim N^*$ according
to \cite{bbfk} 1.7, we see by descending
induction on $i$ that the reduction of $F^{\ge
i} \otimes E_{\sigma (0)}$ has weights at most $
(s+k-1)+t = n+k-1$. The case $i = 0$ gives the
claim.
\end{proof}

This finishes the proof of both, \ref{hltdeform}
and \ref{HRRgerm}.
\end{proof}

{\small

}

\begin{flushleft}
\textsc{G. Barthel, L. Kaup \\ Fachbereich Mathematik und Statistik \\
Universit\"at Konstanz Fach D 203 \\
D-78457 Konstanz \\ Germany}\\
\textit{E-mail adresses}: {\tt
Gottfried.Barthel@uni-konstanz.de,
Ludger.Kaup@uni-konstanz.de}
\end{flushleft}

\begin{flushleft}
\textsc{J.P. Brasselet \\ IML/CNRS, Luminy Case
907 \\
F-13288 Marseille Cedex 9 \\ France}\\
\textit{E-mail adress}: {\tt
jpb@iml.univ-mrs.fr}

\end{flushleft}

\begin{flushleft}
\textsc{K.-H. Fieseler \\ Matematiska Institutionen Box 480 \\
Uppsala Universitet \\SE-75106 Uppsala \\
Sweden}\\
\textit{E-mail adress}: {\tt khf@math.uu.se}

\end{flushleft}


\begin{thebibliography}{99}
\bibitem[BBFK$_{1}$]{hirz} {\sc G.~Barthel, J.-P.~Brasselet,
K.-H.~Fieseler and L.~Kaup,} {\it Equivariant Intersection
Cohomology of Toric Varieties}, Algebraic Geometry, Hirzebruch 70,
45--68, Contemp. Math. {\bf 241}, Amer. Math. Soc., Providence,
R.I., 1999.

\bibitem[BBFK$_{2}$]{bbfk} {\sc G.~Barthel, J.-P.~Brasselet,
K.-H.~Fieseler and L.~Kaup,} {\it Combinatorial Intersection
Cohomology for Fans}, T\^ohoku Math.~J. {\bf 54} (2002), 1--41.

\bibitem[BBFK$_{3}$]{mars} {\sc G.~Barthel, J.-P.~Brasselet,
K.-H.~Fieseler and L.~Kaup,} {\it Combinatorial Duality and
Intersection Product: A Direct Approach}, (pr)e-print {\tt
math.AG/0309352v1} (21~pages), 2003, to appear in  T\^ohoku
Math.~J.

\bibitem[BreLu$_1$]{brel1} {\sc P.~Bressler and V.~Lunts,} {\it
Intersection cohomology on nonrational
polytopes}, Compos. Math. {\bf 135} (2003),
245--278.

\bibitem[BreLu$_2$]{brel2} {\sc P.~Bressler and V.~Lunts,} {\it Hard
Lefschetz theorem and Hodge-Riemann relations for intersection
cohomology of nonrational polytopes}, (pr)e-print {\tt
math.AG/0302236~v2} (46~pages), 2003.

\bibitem[Bri]{bri} {\sc M.~Brion,} {\it The structure of the polytope
algebra}, T\^ohoku Math.~J. {\bf 49} (1997),
1--32.

\bibitem[Ka]{karu} {\sc K.~Karu,} {\it Hard Lefschetz Theorem for
Nonrational Polytopes}, Invent. Math. {\bf 157}
(2004), 419--447.

\bibitem[Mc]{mcmu} {\sc P.~McMullen,} {\it On simple Polytopes},
Invent.~Math. {\bf 113} (1993), 419--444.

\bibitem[St]{stan} {\sc R.~Stanley,} {\it Generalized h-vectors,
intersection cohomology of toric varieties and related results},
{\sc M.~Nagata, H.~Matsumura}, eds., {\it Commutative Algebra and
Combinatorics}, 187--213, Adv.\ Stud.\ Pure Math.~{\bf 11},
Kinokunia, Tokyo, and North Holland, Amsterdam/New York, 1987.

\bibitem[Ti]{timo} {\sc V. A.~Timorin,} {\it An analogue of the
Hodge-Riemann relations for simple convex polytopes}, Russian
Mathematical Surveys {\bf 54.2} (1999), 381--426

\end{thebibliography}
\end{document}